\documentclass[12pt]{article}
\usepackage{amsfonts}
\usepackage{latexsym,amssymb}
\usepackage{mathrsfs}
\usepackage{graphicx}
\usepackage{psfrag}
\usepackage{amsmath, amsbsy}
\usepackage{amsopn, amstext}
\usepackage{multicol}
\usepackage{latexsym,amssymb}

\allowdisplaybreaks

\usepackage{hyperref}

\def\sqr#1#2{{\vcenter{\vbox{\hrule height.#2pt
              \hbox{\vrule width.#2pt height#1pt \kern#1pt \vrule width.#2pt}
              \hrule height.#2pt}}}}
\def\signed #1{{\unskip\nobreak\hfil\penalty50
              \hskip2em\hbox{}\nobreak\hfil#1
              \parfillskip=0pt \finalhyphendemerits=0 \par}}
\def\endpf{\signed {$\sqr69$}}

\def\dbR{{\mathop{\rm l\negthinspace R}}}

\def\3n{\negthinspace \negthinspace \negthinspace }
\def\2n{\negthinspace \negthinspace }
\def\1n{\negthinspace }

\def\ds{\displaystyle}

\def\dbR{{\mathbb{R}}}

\def\={\buildrel \triangle \over =}

\def\mE{\mathbb{E}}

\def\b{\beta}

\def\d{\delta}

\def\l{\lambda}

 \def\n{\nabla}
\def\si{\sigma}
\def\t{\times}
\def\f{\varphi}
\def\th{\theta}
\def\o{\omega}

\def\ov{\bar v}
\def\ns{\noalign{\ss} }

\def\pa{\partial}
\def\G{\Gamma}
\def\D{\Delta}

\def\Si{\Sigma}

\def\O{\Omega}

\def\cF{{\cal F}}

\def\cL{{\cal L}}

\def\cP{{\cal P}}

\def\cX{{\cal X}}

\def\mE{{\mathbb{E}}}

\def\no{\noindent}

\def\ms{\medskip}
\def\bs{\bigskip}
\def\q{\quad}
\def\qq{\qquad}
\def\hb{\hbox}

\def\max{\mathop{\rm max}}
\def\min{\mathop{\rm min}}

\def\pa{\partial}

\def\cd{\cdot}

\def\div{\hbox{\rm div$\,$}}

\def\|{\Big |}
\def\({\Big (}
\def\){\Big )}
\def\[{\Big[}
\def\]{\Big]}
\def\be{\begin{equation}}
\def\bel{\begin{equation}\label}
\def\ee{\end{equation}}
\def\bt{\begin{theorem}}
\def\bcd{\begin{condition}}
\def\ecd{\end{condition}}
\def\et{\end{theorem}}
\def\bc{\begin{corollary}}
\def\ec{\end{corollary}}
\def\bde{\begin{definition}}
\def\ede{\end{definition}}
\def\bl{\begin{lemma}}
\def\el{\end{lemma}}
\def\bp{\begin{proposition}}
\def\ep{\end{proposition}}
\def\br{\begin{remark}}
\def\er{\end{remark}}
\def\ba{\begin{array}}
\def\ea{\end{array}}
\def\ed{\end{document}}
\def\ns{\noalign{\ms}}
\def\ds{\displaystyle}

\def\square#1{\vbox{\hrule\hbox{\vrule height#1%
     \kern#1\vrule}\hrule}}
\def\rectangle#1#2{\vbox{\hrule\hbox{\vrule height#1%
     \kern#2\vrule}\hrule}}

\font\tenbb=msbm10 \font\sevenbb=msbm7 \font\fivebb=msbm5

\newfam\bbfam
\scriptscriptfont\bbfam=\fivebb \textfont\bbfam=\tenbb
\scriptfont\bbfam=\sevenbb

\newtheorem{lemma}{Lemma}[section]
\newtheorem{remark}{Remark}[section]

\newtheorem{theorem}{Theorem}[section]
\newtheorem{corollary}{Corollary}[section]

\newtheorem{definition}{Definition}[section]
\newtheorem{proposition}{Proposition}[section]
\newtheorem{condition}{Condition}[section]

\makeatletter
   
   \@addtoreset{equation}{section}
\makeatother

\begin{document}
\title{\bf Exact Controllability for  Stochastic Schr\"{o}dinger Equations\thanks{This work
is partially supported by the NSF of China
under grant 11101070 and the Fundamental
Research Funds for the Central Universities in
China under grants ZYGX2012J115.  \ms} \ms}

\author{Qi L\"u\thanks{Universit\'{e} Pierre et Marie Curie-Paris VI
UMR 7598, Laboratoire Jacques-Louis Lions, 4, place Jussieu
Paris, F-75005 France; School of Mathematical Sciences, University of Electronic
 Science and Technology of China, Chengdu, 610054, China.
  {\small\it E-mail:} {\small\tt luqi59@163.com}.} }

\date{}

\maketitle

\begin{abstract}\no
This paper is addressed to studying the
exact controllability for  stochastic
Schr\"{o}dinger equations  by two
controls. One is a boundary control in
the drift term and the other is an
internal control in the diffusion term.
By means of the standard duality
argument, the control problem is
converted into  an observability
problem for backward stochastic
Schr\"{o}dinger equations, and the desired
observability estimate is obtained by a
global Carleman estimate. At last, we
give a result  about the lack of exact
controllability, which shows that the action
of two controls is necessary.
\end{abstract}

\bs

\no{\bf 2010 Mathematics Subject
Classification}.  Primary 93B05; Secondary
93B07, 93E20, 60H15.

\bs

\no{\bf Key Words}. Stochastic
Schr\"{o}dinger equations, exact
controllability, observability,
Carleman estimate.

\ms

\section{Introduction}

\q Let $T > 0$, $G \subset
\mathbb{R}^{n}$ ($n \in \mathbb{N}$) be
a given bounded domain with the $C^{2}$
boundary $\G$. Let $\G_0$ be a suitable
chosen nonempty subset (to be given
later) of $\G$. Put $Q \= (0,T) \t G$,
$\Si \= (0,T) \t \G$  and $\Si_0 \=
(0,T) \t \G_0$.

 Let
$(\O,\cF,\{\cF_t\}_{t\ge 0},P)$ be a
complete filtered probability space on
which a one dimensional standard
Brownian motion $\{B(t)\}_{t\ge0}$ is
defined  such that $\{\cF_t\}_{t\ge 0}$
is the natural filtration generated by
$\{B(t)\}_{t\ge0}$, augmented by all
the $P$-null sets in $\cF$. Let $H$ be
a Banach space. We denote by
$L_{\cF}^2(0,T;H)$ the Banach space
consisting of all $H$-valued
$\{\cF_t\}_{t\ge 0}$-adapted processes
$X(\cd)$ such that
$\mathbb{E}(|X(\cd)|_{L^2(0,T;H)}^2)<\infty$;
by $L_{\cF}^\infty(0,T;H)$ the Banach
space consisting of all $H$-valued
$\{\cF_t\}_{t\ge 0}$-adapted bounded
processes; and by
$C_{\cF}([0,T];L^2(\O;H))$ the Banach
space consisting of all $H$-valued
$\{\cF_t\}_{t\geq 0}$-adapted
processes $X(\cd)$ such that $|X(\cd)|_{L^2(\O;H)}\in C([0,T])$. All of the above spaces
are endowed  with the canonical norm.

Denote by $\nu(x)$ the unit outward
normal vector of $G$ at $x\in \G$.  Let
$x_0\in\big(\mathbb{R}^n\setminus
\overline G\big)$. In what follows, we
choose
\begin{equation}\label{G0}
\G_0=\big\{ x\in \G :\, (x-x_0)\cdot
\nu(x)>0  \big\}.
\end{equation}

The main purpose of this paper is to
study the exact controllability of the
following controlled linear stochastic
Schr\"{o}dinger equation
\begin{equation}\label{system1}
\left\{
\begin{array}{lll}
\ds idy + \D ydt = (a_1 \cdot \nabla y + a_2 y + f)dt + (a_3 y + g)dB(t)  &\mbox { in } Q, \\
\ns\ds y = 0 &\mbox{ on } \Si\setminus\Si_0, \\
\ns\ds y=u &\mbox{ on } \Si_0, \\
y(0) = y_0 &\mbox{ in } G.
\end{array}
\right.
\end{equation}
Here, the initial state $y_0 \in
L^2(\O,\cF_0,P;H^{-1}(G))$, the control
$$
u\in L^2_\cF(0,T;L^2(\G_0)),\q g\in
L^2(0,T;H^{-1}(G)),$$ the
nonhomogeneous term $f\in
L^2_\cF(0,T;L^2(G))$ and the
coefficients $a_i$($i=1,2,3$) satisfy
\begin{equation}\label{ai}
\left\{\begin{array} {ll} \ds  ia_1 \in
L_{
\mathcal{F}}^{\infty}(0,T;W^{2,\infty}(G;\mathbb{R}^{n})\cap
W_0^{1,\infty}(G;\mathbb{R}^{n})),
\\
\ns\ds a_2
\in L_{ \mathcal{F}}^{\infty}(0,T;W^{1,\infty}(G)), \\
\ns  \ds  a_3 \in L_{\mathcal{
F}}^{\infty}(0,T;W^{1,\infty}(G)).
\end{array}
\right.
\end{equation}

System \eqref{system1} has a
nonhomogeneous boundary condition. As
the deterministic nonhomogeneous
boundary problem,  the solution to
\eqref{system1} is understood in the
transposition sense. Hence,  we first
introduce the following backward
stochastic Schr\"{o}dinger equation
\begin{equation}\label{system2}
\left\{
\begin{array}{lll}
\ds idz + \D zdt = \big(b_1 \cdot
\nabla z
+ b_2 z + b_3Z \big)dt + ZdB(t)  &\mbox { in } (0,\tau)\times G, \\
\ns\ds z = 0 &\mbox{ on } (0,\tau)\times \G, \\
z(\tau) = z_\tau &\mbox{ in } G,
\end{array}
\right.
\end{equation}
where $z_T\in L^2_{\cF_T}(\O;H_0^1(G))$,
the coefficients $b_i$($i=1,2,3$)
satisfy
\begin{equation}\label{bi}
\left\{\begin{array} {ll} \ds  ib_1 \in
L_{
\mathcal{F}}^{\infty}(0,T;W_0^{1,\infty}(G;\mathbb{R}^{n})),
\\
\ns\ds b_2
\in L_{\mathcal{F}}^{\infty}(0,T;W^{1,\infty}(G)), \\
\ns  \ds  b_3 \in L_{\mathcal{
F}}^{\infty}(0,T;W^{1,\infty}(G)).
\end{array}
\right.
\end{equation}
For the convenience of the reader, we
recall the definition of the solution
to \eqref{system2} first.

\begin{definition}\label{def bt sol}
A solution to the equation
\eqref{system2} is a pair of stochastic
processes
$$(z,Z)\in L^\infty_{\cF}(0,T;H^1_0(G))\t
L^2_{\cF}(0,T;H^1_0(G))$$ such that for
every $\psi\in C_0^\infty(G)$ and a.e.
$(t,\o)\in [0,\tau]\times\O$, it holds
that
\begin{equation}\label{def id1}
\begin{array}{ll}
\ds \q\int_G  iz_\tau(x)\psi(x)dx -
\int_G iz(t)\psi(x)dx  -
\int_t^\tau\int_G \nabla
z(s,x)\cd\nabla\psi(x)dxds
\\ \ns\ds = \!\int_t^\tau \!\!\int_G
\big[b_1(s,x)\nabla z(s,x)
 + b_2(s,x)z(s,x)+ b_3(s,x)Z(s,x)\big]\psi(x)dxds \\
\ns\ds \q + \int_t^\tau\int_G  Z(s,x)
\psi(x)dxdB(s).
\end{array}
\end{equation}
\end{definition}

\medskip

Let us recall the well-posedness result
of \eqref{system2}(see
\cite{Al-Hussein1,Mahmudov1} for the
proof).

\begin{lemma}\label{well posed1}
For any $z_\tau\in
L^2(\O,\cF_\tau,P;H_0^1(G))$, the
equation \eqref{system2} admits a
unique solution $(z,Z)$. Moreover,
$(z,Z)$ satisfies that
\begin{equation}\label{best1}
|z|_{L^\infty_\cF(0,\tau;H_0^1(G))} +
|Z|_{L^2_\cF(0,\tau;H^1_0(G))} \leq
e^{Cr_1}
|z_\tau|_{L^2(\O,\cF_\tau,P;H_0^1(G))},
\end{equation}
where
$$ r_1\=
|b_1|^2_{L^\infty_{\cF}(0,T;W^{1,\infty}(G;\dbR^n))}
+
\sum_{i=2}^3|b_i|^2_{L^\infty_{\cF}(0,T;W^{1,\infty}(G))}
+1.
$$
\end{lemma}
Here and in the sequel, we will use $C$
to denote a generic positive constant depends on $T$, $G$, $\G_0$(unless otherwise stated),
which may vary from line to line.

Further, we need the following  result for the hidden
regularity of the solution to
\eqref{system2}.

\begin{proposition}\label{prop-hid}
Let $z$ be a solution to
\eqref{system2}, then $\frac{\pa
z}{\pa\nu}|_{\G}\in
L^2_{\cF}(0,\tau;L^2(\G))$. Further, we
have the following estimate
\begin{equation}\label{hid-eq1}
\Big| \frac{\pa z}{\pa\nu}
\Big|_{L^2_\cF(0,\tau;L^2(\G))} \leq
e^{Cr_1}|z_\tau|_{L^2_{\cF_\tau}(\O;H_0^1(G))}.
\end{equation}
\end{proposition}

\begin{remark}
Proposition \ref{prop-hid} shows that,
solutions of \eqref{system2} enjoy a
higher regularity on the boundary than
the one provided by the classical trace
theorem of Sobolev spaces. Such kind
of result  is called hidden regularity
of the solution. There are a great many
studies in this topic for
deterministic partial differential
equations in the literature(see
\cite{La,Lions1} for example).
\end{remark}

Now we can give the definition of the
 solution  to
\eqref{system1}.

\begin{definition}\label{def ft sol}
A solution to the system
\eqref{system1} is a process $y\in
C_{\cF}([0,T];L^2(\O;H^{-1}(G)))$ such
that for every $\tau\in [0,T]$ and
every $z_\tau\in
L^2_{\cF_\tau}(\O;H_0^1(G))$ it holds
that
\begin{equation}\label{def id}
\begin{array}{ll}
\ds \q\mE\int_G y(\tau,x)z_\tau(x)dx -
\mE\int_G y_0(x)z(0)dx
\\ \ns\ds = \mE\int_0^\tau \int_{\G_0}\frac{\pa z}{\pa\nu}udxds
+ \mE\int_0^\tau\int_G fzdxdt +
\mE\int_0^\tau\int_G gZdxdt.
\end{array}
\end{equation}
Here $(z,Z)$ solves \eqref{system2}
with
$$
b_1 = -a_1,\q b_2 = -\div(a_1)+a_2,
b_3=-a_3.
$$
\end{definition}

\begin{remark}
The solution to \eqref{system1} is
defined in the transposition sense. It
is well studied that such kind of
solutions for deterministic
nonhomogeneous boundary value problems
in the literature(see
\cite{Lions1,Lions2} for example). On
the other hand, the stochastic
counterpart is almost open. We only
consider a very special problem in this
topic. The study of general stochastic
nonhomogeneous boundary value problems
are very interesting but difficult
problems, which is far beyond the scope
of this paper.
\end{remark}

We have the following well-posedness
result for \eqref{system1}.

\begin{proposition}\label{well posed}
For each $y_0\in
L^2_{\cF_0}(\O;H^{-1}(G))$, the system
\eqref{system1} admits a unique
solution $y$. Further,   for every
$y_0\in L^2_{\cF_0}(\O;H^{-1}(G))$, it
holds that
\begin{equation}\label{well posed est}
\begin{array}{ll}\ds
\q|y|_{C_{\cF}([0,T];L^2(\O;H^{-1}(G)))}\\
\ns\ds \leq e^{Cr_1}\big(
\mE|y_0|_{H^{-1}(G)} +
|f|_{L^2_{\cF}(0,T;L^2(G))} +
|u|_{L^2_{\cF}(0,T;L^2(\G_{0}))} +
|g|_{L^2_{\cF}(0,T;H^{-1}(G))}\big).
\end{array}
\end{equation}
Here $$
r_2 \ = |a_1|^2_{L_{
\mathcal{F}}^{\infty}(0,T;W_0^{1,\infty}(G;\mathbb{R}^{n}))}
 +
|a_2|^2_{L_{\mathcal{F}}^{\infty}(0,T;W^{1,\infty}(G))}
 +
|a_3|^2_{L_{\mathcal{
F}}^{\infty}(0,T;W^{1,\infty}(G))}  +
1.$$
\end{proposition}

\medskip

Now we can give the definition of the
exact controllability of
\eqref{system1}.

\begin{definition}\label{exact def}
System \eqref{system1} is said to be
exactly controllable at time $T$ if for
every initial state $y_0\in
L^2_{\cF_0}(\O;H^{-1}(G))$ and every
$y_1\in L^2_{\cF_T}(\O;H^{-1}(G))$, one
can find a pair of controls $(u,g)\in
L^2_\cF(0,T;L^2(\G_{0}))\t
L^2_{\cF}(0,T;H^{-1}(G))$ such that the
solution $y$ of the system
\eqref{system1} satisfies that $y(T) =
y_1$ in $L^2_{\cF_T}(\O;H^{-1}(G))$.
\end{definition}

As the deterministic case,  the exact
controllability of \eqref{system1} can
be reduced to an observability estimate
of its dual system, that is, the
equation \eqref{system2} with $\tau=T$ in our case.
For the latter one, we will prove the
following result.

\begin{theorem}\label{observability}
All solutions of the equation
\eqref{system2} satisfy that
\begin{equation}\label{exact ob est}
|z_T|_{L^2_{\cF_T}(\O;H^1_0(G))}\leq
e^{Cr_1}\big(|z|_{L_\cF^2(0,T;L^2(\G_0))}
+ |Z|_{L^2_\cF(0,T;H^1_0(G))}\big).
\end{equation}
\end{theorem}

By means of Theorem
\ref{observability},  we can obtain the
following the exact
controllability result of the system
\eqref{system1}.

\begin{theorem}\label{exact th}
System \eqref{system1} is exactly
controllable at any time $T>0$.
\end{theorem}

Further, we also have the following
result about the lack of exact
controllability result for
\eqref{system1} if the control $u$ in
the drift term is zero.

\begin{theorem}\label{non control th}
If $u\equiv 0$  in  \eqref{system1},
then  \eqref{system1} is not exactly
controllable at any time $T$.
\end{theorem}

\begin{remark}
According to the lack of exact
controllability for linear stochastic
ordinary differential equations in
\cite{Peng}, it is clear that the
internal control $g$ is necessary for
the exact controllability of the system
\eqref{system1}. From Theorem \ref{non
control th}, we know that $u$ is also
necessary. Hence, we have that one
should utilize both $u$ and $g$ to get
the exact controllability of
\eqref{system1}.
\end{remark}

There exist many approaches and results
addressing the controllability problem
for determinisitc Schr\"{o}dinger
equations.  \cite{Zuazua} is a nice
survey for the works in this respect
before 2002. For the works after 2002,
we refer the readers to
\cite{BL,BBZ,BZ,Cazacu,EP,RZ} and the
rich references therein.   However,
people know very little about the
stochastic counterpart. To our best
knowledge, there is no published result
for this problem.

Generally speaking, there are four main
methods for the exact controllability
of  deterministic Schr\"{o}dinger
equations.

The first one is the classical
Rellich-type multiplier approach
(\cite{Machtyngier}). It can be applied
to treat Schr\"{o}dinger equations with
no lower order terms or lower order
terms with constant or small
coefficients. On the other hand, it seems
that it cannot be used to solve our
problem since we do not assume that the
coefficients of lower order terms  are constant or small.

 The second one
is the microlocal analysis approach
(\cite{Lebeau}). This method was first
introduced in
\cite{Bardos-Lebeau-Rauch1} for
obtaining the exact controllability of
wave equations. It is useful to solve
the exact controllability problem for
many kinds of partial differential equations
such as wave equations, Schr\"{o}dinger
equations and plate equations. However,
it seems that there are lots of
obstacles needing to be surmounted if one utilize it to study the stochastic control problem(see
remarks in Section \ref{comment} for
more details).

The third one is based on the Ingham
type inequality(\cite{Komornik1}). This
method works well for Schr\"{o}dinger
equations involved in some special
domains, i.e., intervals, rectangles and
balls. However, it seems that it is
very hard to applied to equations in
general domains.

The last one is the global Carleman
estimate(\cite{Lasiecka-Triggiani-Zhang}).
It can be regarded as a more developed
version of the classical multiplier
method. With respect to the method of
multipliers, the Carleman approach has
the advantage of being more flexible
and allowing to address variable
coefficients,  with respect to the
microlocal one, that it requires less
regularity on coefficients and domain,
and respect to the method based on
Ingham type inequality, that is
requires less restrictions to the
domain. Further, it is robust with
respect to the lower order terms and
can be used to get explicit bounds on
the observability constant/control cost
in terms of the potentials entering in
it.  This is particularly important
when dealing with nonlinear problems by
combing linearization and
fixed point techniques.

Similar to the deterministic setting,
we  use a stochastic version of the
global Carleman estimate to derive
Theorem \ref{observability}. For
this, we borrow some idea from the
proof of the observability estimate for
deterministic Schr\"{o}dinger equations
(see \cite{Lasiecka-Triggiani-Zhang}
for example). However, the stochastic
setting will produce some more troubles. We
cannot simply mimic the method in
\cite{Lasiecka-Triggiani-Zhang} to
solve our problem. To handle these troubles, we choose a different weight function for the Carleman estimate.

\section{Some preliminaries}

\q\,\, In this section, we give some
preliminary results.

First, we prove the well-posedness of
\eqref{system1} in the sense of
Definition \ref{def ft sol}.

\vspace{0.2cm}

{\it Proof of Proposition \ref{well
posed}}\,: {\bf Uniqueness of the
solution.}\; Suppose there are
$y_1(\cd)$ and  $y_2(\cd)$ belong to
$C_\cF([0,T];L^2(\O;H^{-1}(G)))$ such
that \eqref{def id} holds. Then, we see
$$
\mE\int_G y_1(\tau,x)z_\tau(x)dx =
\mE\int_G y_2(\tau,x)z_\tau(x)dx
\q\mbox{ for all } z_\tau\in
L^2_{\cF_\tau}(\O;H_0^1(G)).
$$
This concludes that $y_1=y_2$.

\vspace{0.2cm}

{\bf Existence of the solution.} \; Let
us define a linear functional $F$ on
$L^2_{\cF_\tau}(\O;H_0^1(G))$ as
$$
F(z_\tau) =
(y_0,z(0))_{H^{-1}(G),H_0^1(G)} +
\mE\int_0^\tau \int_{\G_0}\frac{\pa
z}{\pa\nu}udxds + \mE\int_0^\tau\int_G
fzdxdt + \mE\int_0^\tau\int_G gZdxdt.
$$
It is easy to show that $F$ is a
bounded linear functional on
$L^2_{\cF_\tau}(\O;H_0^1(G))$. By Riesz
Representation Theorem, we know there
exists a $y_\tau \in
L^2_{\cF_\tau}(\O;H^{-1}(G))$ such that
$$
F(z_\tau) =
\mE(y_\tau,z_\tau)_{H^{-1}(G),H_0^1(G)}.
$$
Define a process $y(\cd)$ by
$y(\tau)=y_\tau$. Now we prove that
$y(\cd)\in
C_\cF([0,T];L^2(\O;H^{-1}(G)))$.

\vspace{0.1cm}

Let $\xi\in L^2_{\cF_T}(\O;H_0^1(G))$.
Consider the following forward random
Schr\"{o}dinger equation
\begin{equation}\label{system3}
\left\{
\begin{array}{lll}
\ds id\tilde z + \D \tilde zdt = \big(-a_1 \cdot \nabla \tilde z-\div a_1 \tilde z + a_2 \tilde z \big)dt   &\mbox { in } (\tau,\tau+\d)\times G, \\
\ns\ds \tilde z = 0 &\mbox{ on } (\tau,\tau+\d)\times \G, \\
\tilde z(\tau) = \mE(\xi|\cF_{\tau})
&\mbox{ in } G,
\end{array}
\right.
\end{equation}
It is easy to see that
\begin{equation}\label{4.16-eq1}
\lim_{\d\to 0^+} \mE|\tilde
z(\tau+\d)-\tilde
z(\tau)|^2_{H_0^1(G)}=0.
\end{equation}
Further, since $\{\cF_t\}_{t\geq
0}$ is the
natural filtration of $\{B(t)\}_{t\geq
0}$, we have
\begin{equation}\label{4.16-eq2}
\lim_{\d\to 0^+}
\mE\big|\mE(\xi|\cF_{\tau+\d}) -
\mE(\xi|\cF_{\tau})\big|^2_{H_0^1(G)}=0.
\end{equation}
From \eqref{4.16-eq1} and
\eqref{4.16-eq2}, we see
\begin{equation}\label{4.16-eq3}
\lim_{\d\to 0^+} \mE|\tilde
z(\tau+\d)-\mE(z|\cF_{\tau+\d})|^2_{H_0^1(G)}=0.
\end{equation}
Let $(z_1(\cd),Z_1(\cd))$ and
$(z_2(\cd),Z_2(\cd))$ be the solution
to  \eqref{system2} with final data
$\mE(\xi|\cF_{\tau+\d})$ and $\tilde
z(\tau+\d)$, respectively. Then, from
Lemma \ref{well posed1} and
Proposition \ref{prop-hid}, we know
\begin{equation}\label{3.25-eq1}
\left\{
\begin{array}{ll}\ds
\lim_{\d\to 0^+} \big| z_1-
z_2\big|_{L^2_\cF(0,\tau;L^2(G))}=0,\\
\ns\ds \lim_{\d\to 0^+}
\Big|Z_1-Z_2\Big|_{L^2_\cF(0,\tau;L^2(\G_0))}=0,\\
\ns\ds \lim_{\d\to 0^+} \Big|\frac{\pa
z_1}{\pa\nu}-\frac{\pa
z_2}{\pa\nu}\Big|_{L^2_\cF(0,\tau;L^2(\G_0))}=0,\\
\ns\ds \lim_{\d\to 0^+}
|z_1(0)-z_2(0)|_{H_0^1(G)}=0.
\end{array}
\right.
\end{equation}
Denote by $(z_3,Z_3)$ the solution to
\eqref{system2} with the final datum
$z_3(\tau)=\mE(\xi|\cF_{\tau})$. From
the uniqueness of the solution to
\eqref{system3} and \eqref{system2}, we
find that
\begin{equation}\label{3.25-eq2}
z_3 = z_2 \mbox{ in } [0,\tau]\times G
\mbox{ and } Z_3 = Z_2 \mbox{ in }
[0,\tau]\times G.
\end{equation}
From the definition of the solution to
\eqref{system1}, we have
$$
\begin{array}{ll}\ds
\q\mE\big(y(\tau+\d)-y(\tau),\xi\big)_{H^{-1}(G),H_0^1(G)}\\
\ns\ds = \mE\big(y(\tau+\d),\xi\big)_{H^{-1}(G),H_0^1(G)} - \mE\big(y(\tau),\xi\big)_{H^{-1}(G),H_0^1(G)}\\
\ns\ds =
\mE\big(y(\tau+\d),\mE(\xi|\cF_{\tau+\d})\big)_{H^{-1}(G),H_0^1(G)}
- \mE\big(y(\tau),\mE(\xi|\cF_{\tau})\big)_{H^{-1}(G),H_0^1(G)}\\
\ns\ds =
\mE\big(y_0,z_1(0)-z_3(0)\big)_{H^{-1}(G),H_0^1(G)}
+ \mE\int_0^\tau \int_{\G_0}\(\frac{\pa
z_1}{\pa\nu} - \frac{\pa
z_3}{\pa\nu}\)ud\G ds\\
\ns\ds \q  +\mE\int_0^\tau\int_G
f(z_1-z_3)dxdt + \mE\int_0^\tau\int_G
g(Z_1-Z_3)dxdt + \mE\int_\tau^{\tau+\d}
\int_{\G_0} \frac{\pa z_1}{\pa\nu} ud\G
ds\\
\ns\ds \q +\mE\int_\tau^{\tau+\d}\int_G
fz_1 dxdt +
\mE\int_\tau^{\tau+\d}\int_G g Z_1
dxdt.
\end{array}
$$
This, together with \eqref{3.25-eq1}
and \eqref{3.25-eq2}, implies that
$$
\lim_{\d\to
0^+}\mE\big(y(\tau+\d)-y(\tau),\xi\big)_{H^{-1}(G),H_0^1(G)}=0,\q\mbox{
for any } \xi\in
L^2_{\cF_T}(\O;H_0^1(G)).
$$
Similarly, we can show that
$$
\lim_{\d\to
0^-}\mE(y(\tau+\d)-y(\tau),\xi)_{H^{-1}(G),H_0^1(G)}=0,\q\mbox{
for any } \xi\in
L^2_{\cF_T}(\O;H_0^1(G)).
$$
Hence, we see $y(\cd)\in
C_\cF([0,T];L^2(\O;H^{-1}(G)))$.
\endpf

\vspace{0.2cm}

Next, for the sake of completeness, we
give an energy estimate for the
equation \eqref{system2}.

\medskip
\begin{proposition} \label{prop1}
For all $z$ which solve the equation
\eqref{system2}, it holds that
\begin{equation}\label{energyesti1}
\mathbb{E}|z(t)|^2_{H_0^1(G)} \leq e^{C
r_1} \Big(\mathbb{E}|z(s)|^2_{H^1_0(G)}
+ |Z|^2_{L^2_{\cF}(0,T;H_0^1(G))}\Big),
\end{equation}
for any $0\leq s \leq t \leq\tau$.
\end{proposition}
\medskip

{\it Proof }: By direct computation, we
have
\begin{equation}\label{Eyt}
\begin{array}{ll} \ds \q\mathbb{E}|z(t)|^2_{ L^2(G)} - \mathbb{E}|z(s)|^2_{
L^2(G)}\\
\ns  \ds =  \mathbb{E}\int_s^t\int_G (z
d\bar{z}+\bar{z}dz + dz
d\bar{z})dx\\
\ns  \ds =  \mathbb{E}\int_s^t\int_G
\Big\{-i z\big(\D \bar{z} + b_1\cdot
\nabla\bar{z}  - b_2
\bar{z} - b_3\overline Z \big)\\
\ns\ds \qq\qq\q + i\bar{z}\big(\D z -
b_1\cdot
\nabla z  - b_2 z - b_3 Z  \big)+ Z\overline Z \Big\}dxd\si \\
\ns  \ds   =\mathbb{E}\int_s^t\int_G
\Big\{-i[\div(z\nabla\bar{z})-|\nabla
z|^2 + \div(|z|^2 b_1) -\div(b_1)|z|^2  - b_2|z|^2
- b_3 z\overline Z ] \\
\ns  \ds \qq \qq + i
[\div(\bar{z}\nabla z)-|\nabla z|^2  - b_2|z|^2- b_3\bar z Z ] + |Z|^2 \Big\}dxd\si\\
\ns  \ds \leq\mathbb{E}\int_s^t
2\Big[(|b_1|_{W^{1,\infty}(G;\dbR^n)}+|b_3|_{L^{\infty}(G)}+1)|z|^2_{L^2(G)}+
|Z|_{L^2(G)}^2 \Big]dxd\si
\end{array}
\end{equation}
and
\begin{equation}\label{Etyt}
\begin{array}{ll}
\ds\q\mathbb{E}|\nabla z(t)|^2_{
L^2(G)} - \mathbb{E}|\nabla
z(s)|^2_{L^2(G)}\\
\ns \ds  =   \mathbb{E}\int_s^t\int_G
(\nabla z d\n\bar{z} + \nabla
\bar{z} d\n z + d\nabla z d\nabla \bar{z})dx\\
\ns  \ds =  \mathbb{E}\int_s^t\int_G
\Big\{ \div(\nabla z d\bar{z}) - \D z
d\bar{z} + \div(\nabla\bar{z} dz) - \D
\bar{z}dz +
d\nabla z d\nabla \bar{z} \Big\}dx\\
\ns \ds  =  \mathbb{E}\int_s^t\int_G
\Big\{\D z \Big[i(\D \bar{z} + b_1\cdot
\nabla\bar{z}  - b_2 \bar{z} -
b_3\overline Z )\Big]
\\
\ns\ds \qq\qq\q\; -\D\bar{z}\Big[
i\big(\D z - b_1\cdot \nabla z  - b_2 z
- b_3 Z  \big)\Big]  +|\nabla Z|^2
\Big\}dxd\si \\
\ns  \ds \leq  2\mathbb{E}\int_t^s
\Big\{\big(|b_1|^2_{W^{1,\infty}(G;\mathbb{R}^n)}+|b_3|^2_{W^{1,\infty}(G)}+1\big)|\nabla
z|^2_{L^2(G)}\\
\ns \ds  \qq \qq\; +
\big(|b_2|^2_{W^{1,\infty}(G)}+|b_3|^2_{W^{1,\infty}(G)}+1\big)|z|^2_{L^2(G)}
+|Z|^2_{H_0^1(G)}\Big\} d\si.
\end{array}
\end{equation}
From \eqref{Eyt} and \eqref{Etyt}, we
get
\begin{equation}\label{energyest2}
\begin{array}{ll}\ds
 \q\mathbb{E}|z(t)|^2_{H_0^1(G)} -
\mathbb{E}|z(s)|^2_{H_0^1(G)} \\
\ns\ds \leq 2(r_1+1)\mathbb{E}\int_t^s
|z(\si)|^2_{H_0^1(G)} d\si +
\mathbb{E}\int_t^s
|Z(\si)|^2_{H_0^1(G)} d\si.
\end{array}
\end{equation}
From \eqref{energyest2}, and thanks to
Gronwall's inequality, we arrive at
\begin{equation}\label{energyest3}
\mathbb{E}| y(t)|^2_{H_0^1(G)}\leq
e^{2(r_1+1)}\Big\{\mathbb{E}|y(s)|^2_{H_0^1(G)}
+ \mathbb{E}\int_0^\tau
|Z|^2_{H_0^1(G)} d\si\Big\},
\end{equation}
which implies the inequality
\eqref{energyesti1} immediately.
\endpf

\begin{remark}
The proof of this proposition is almost
standard.  Indeed, if we regard $z$ as
a solution to a forward stochastic
Schr\"{o}dinger equation with a
nonhomogeneous term $Z$, then it is a
standard energy estimate for such kind
of equation.
\end{remark}

Next, we give a proof of Proposition
\ref{prop-hid}. For this, we first
recall a pointwise identity. For
simplicity, in what follows, we
adopt the notation $\ds z_i \equiv
z_{i}(x) \= \frac{\partial
z(x)}{\partial x_i}$, where $x_i$ is
the $i$-th coordinate of a generic
point $x=(x_1,\cdots, x_n)$ in
$\mathbb{R}^{n}$. In a similar manner,
we  use the notation $y_i$, $v_i$,
etc., for the partial derivatives of
$y$ and $v$ with respect to $x_i$. Let
us recall the following identity.

\medskip

\begin{lemma}\cite[Proposition 2.3]{Lu1}\label{prop2}
Let $\mu = \mu(x) =
(\mu^1,\cdots,\mu^n):\mathbb{R}^n \to
\mathbb{R}^n$ be a vector field of
class $C^1$ and $z$ an
$H^2_{loc}(\mathbb{R}^n)$-valued
$\{\mathcal{F}_t\}_{t\geq 0}$-adapted
process. Then for a.e. $x \in
\mathbb{R}^n$ and P-a.s. $\omega \in
\Omega$, it holds that
\begin{equation}\label{identity2}
\begin{array}{ll}
\q \ds \mu\cdot\nabla\bar{z}(i dz +
\Delta z dt) +
\mu\cdot\nabla z(-i d\bar{z} + \Delta \bar{z} dt)\\
\ns  = \ds  \nabla\cd \Big[
(\mu\cdot\nabla \bar{z})\nabla z+
(\mu\cdot\nabla z)\nabla \bar{z}  - i
(z d\bar{z}) \mu - |\nabla
z|^2\mu \Big]dt + d(i\mu\cd\nabla \bar{z} z)\\
\ns \ds\q - 2\sum_{j,k=1}^n \mu^k_j
z_{j}\bar{z}_{k}dt + (\nabla\cdot \mu)
|\nabla z|^2 dt + i(\nabla\cdot \mu) z
d\bar{z} - i(\mu\cd\nabla d\bar z) dz.
\end{array}
\end{equation}
\end{lemma}

\medskip

By virtue of  Lemma \ref{prop2}, the
proof of Proposition \ref{prop-hid} is
standard. We only give a sketch here.

\vspace{0.2cm}

{\it Sketch of the Proof of Proposition
\ref{prop-hid}} : Since $\Gamma$ is $
C^2$, one can find a vector field
$\mu_0 = (\mu_0^1, \cdots, \mu_0^n) \in
C^1(\overline{G};\mathbb{R}^n)$ such
that $\mu_0 = \nu$ on $\Gamma$(see
\cite[page 18]{Komornik} for the
construction of $\mu_0$). Letting $\mu
= \mu_0$ and $z = y$ in Lemma
\ref{prop2}, integrating it in $Q$ and
taking the expectation, by means of
Proposition \ref{prop2}, with similar
computation in \cite{Zhangxu1},
Proposition \ref{prop-hid} can be
obtained immediately.
\endpf

\vspace{0.2cm}

Next, we recall an identity in the
spirit of \eqref{identity2} but much
more complex, which will play an
important role in establishing the
Carleman estimate for \eqref{system2}.

Let $\b(t,x)\in
C^{2}(\mathbb{R}^{1+n};\mathbb{R})$,
and let $b^{jk}(t,x)\in
C^{1,2}(\mathbb{R}^{1+n};\;\mathbb{R})$
satisfy
\begin{equation}\label{bjk}
b^{jk}=b^{kj},\qq j,k=1,2,\cdots,n.
\end{equation}
Let us define a (formal) second order
stochastic partial differential
operator $\cP$ as
\begin{equation}\label{cp}
\ds
 \cP z \= i\b(t,x)dw+\sum_{j,k=1}^n(b^{jk}(t,x)w_j)_k dt,
 \q i=\sqrt{-1}.
\end{equation}
We have the following equality
concerning $\cP$:

\begin{lemma}\cite[Theorem 3.1]{Lu1}\label{identity1}
Let $\ell,\;\Psi\in
C^2(\mathbb{R}^{1+n};\;\mathbb{R})$ and $\th=e^\ell$.
Assume that $w$ is an
$H^2_{loc}(\mathbb{R}^n,\mathbb{C})$-valued
$\{\cF_t\}_{t\geq 0}$-adapted process.
 Put $ v=\th
w$. Then for a.e. $x
\in \mathbb{R}^n$ and $P$-a.s. $\o\in
\O$,  it holds that

\begin{equation}\label{c2a1}
\begin{array}{ll}
 \ds\th(\cP w\overline {I_1}+\overline{\cP w} I_1)+dM+\div V \\
\ns  =  \ds 2|I_1|^2 dt +\sum_{j,k=1}^n
 c^{jk}(v_k\ov_j+\ov_k v_j)
dt +D|v|^2 dt \\
\ns \ds  +i\sum_{j,k=1}^n\[(\b
b^{jk}\ell_j)_t+
b^{jk}(\b\ell_t)_j\](\ov_kv-v_k\ov) dt
\\
\ns \ds +i\[\b\Psi+\sum_{j,k=1}^n(\b
b^{jk}\ell_j)_k\](\ov
dv-vd\ov)\\
\ns \ds  + (\b^2\ell_t)dvd\ov +
i\sum_{j,k=1}^n \b b^{jk}\ell_j (dv
d\ov_k - dv_kd\ov),
\end{array}
\end{equation}
where
\begin{equation}\label{c2a2}
\left\{
\begin{array}{ll}\ds I_1\= - i\b \ell_t
v - 2\sum_{j,k=1}^n b^{jk}\ell_j v_k +
\Psi v, \\
\ns\ds A\=\sum_{j,k=1}^n
b^{jk}\ell_j\ell_k-\sum_{j,k=1}^n(b^{jk}\ell_j)_k
-\Psi,
\end{array}
\right.
\end{equation}

\begin{equation}
\label{c2a3} \left\{
\begin{array}{ll}\ds
M\=\b^2\ell_t |v|^2 + i\b\sum_{j,k=1}^n b^{jk}\ell_j(\ov_kv-v_k\ov),\\
\ns\ds V\=[V^1,\cdots,V^k,\cdots,V^n],\\
 \ns\ds V^k\=-i \b\sum_{j=1}^n\[b^{jk}\ell_j(vd\ov -\ov
 dv ) + b^{jk}\ell_t(v_j\ov-\ov_jv) dt\]\\
\ns\ds\qq\,\,\,\, - \Psi\sum_{j=1}^n b^{jk}(v_j\ov+\ov_jv) dt + \sum_{j=1}^n b^{jk}(2A\ell_j+\Psi_j)|v|^2 dt \\
\ns\ds\qq\q+\sum_{j,j',k'=1}^n\(2b^{jk'}b^{j'k}-b^{jk}b^{j'k'}\)\ell_j(v_{j'}\ov_{k'}+\ov_{j'}v_{k'})
dt,
\end{array}
\right.
\end{equation}
and
\begin{equation}\label{cc2a3}
\left\{
\begin{array}{ll}\ds
c^{jk}\= \sum_{j',k'=1}^n\[2(b^{j'k}\ell_{j'})_{k'}b^{jk'}-(b^{jk}b^{j'k'}\ell_{j'})_{k'}\] - b^{jk}\Psi,\\
\ns\ds D\=(\b^2\ell_t)_t
+\sum_{j,k=1}^n(b^{jk}\Psi_k)_j+2\[\sum_{j,k=1}^n(b^{jk}\ell_jA)_k+A\Psi\].
\end{array}
\right.
\end{equation}
\end{lemma}

\section{A global Carleman estimate for the equation \eqref{system2}}

In this section, we establish a global
Carleman estimate for the solution to
\eqref{system2}(see Theorem
\ref{thcarleman est} below).

To begin with, let us introduce the
weight functions to be used in our
Carleman estimate. Let
\begin{equation}\label{psi}
\psi(x) = |x-x_0|^2 + \si,
\end{equation}
where $\si$ is a positive constant
such that $\psi \geq
\frac{5}{6}|\psi|_{L^{\infty}(G)}$. Let
$s>0$ and $\l>0$.  Put
\begin{equation}\label{lvarphi}
\ell = s\frac{e^{4\l \psi} - e^{5\l
|\psi|_{L^{\infty}(G)}}}{t^2(T-t)^2},
\qq \varphi = \frac{e^{4\l \psi}
}{t^2(T-t)^2}.
\end{equation}

We have the following  global Carleman
inequality.

\begin{theorem}\label{thcarleman est}
According to  \eqref{lvarphi}, there is
an $s_1>0$ (depending on $r_1$) and a
$\l_1>0$ such that for each $s\geq
s_1$, $\l\geq \l_1$ and for any
solution of the equation
\eqref{system2}, it holds that
\begin{eqnarray}\label{carleman est}
\begin{array}{ll}
\ds \q\mathbb{E}\int_Q
\theta^2\Big(s^3\l^4\varphi^3 |z|^2 +
s\l\varphi
|\nabla z|^2\Big) dxdt \\
\ns \ds \leq  C \[\mathbb{E}\int_Q
\theta^2 \Big(
 s^2\l^2\varphi^2 |Z|^2 + |\nabla Z|^2 \Big)dxdt +
\mathbb{E}\int_0^T\int_{\G_0}\theta^2
s\l\varphi\Big| \frac{\pa z}{\pa
\nu}\Big|^2d\G dt \].
\end{array}
\end{eqnarray}
\end{theorem}

{\it Proof of Theorem {\ref{thcarleman
est}}}: The proof is divided into
three steps.

\medskip

\textbf{Step 1.} We choose $\b = 1$ and
$(b^{jk})_{1\leq j,k\leq n}$ to be the
identity matrix. Put
$$
\d^{jk} =
\left\{
\begin{array}{ll}\ds 1,&\mbox{
if } j=k,\\ \ns\ds 0,&\mbox{ if } j\neq
k.
\end{array}
\right.
$$
Applying Lemma \ref{identity1} to the
equation \eqref{system2} with $\theta$
given by \eqref{lvarphi}, $w$ replaced
by $z$ and $v = \theta z$, we obtain
that
\begin{equation}\label{identity2.1}
\begin{array}{ll}
& \ds\theta\cP z  {( i\b \ell_t \bar{v}
- 2\sum_{j,k=1}^n b^{jk}\ell_j
\bar{v}_k + \Psi \bar{v})} +
\theta\overline{\cP z} {(- i\b \ell_t v
- 2\sum_{j,k=1}^n b^{jk}\ell_j v_k  +
\Psi v)}\\
\ns & \ds \q + \;dM + \div V  \\
\ns = & \ds 2\Big|- i\b \ell_t v -
2\sum_{j,k=1}^n b^{jk}\ell_j v_k  +
\Psi v\Big|^2dt +
\sum_{j,k=1}^nc^{jk}(v_k\ov_j+\ov_k
v_j) dt + D|v|^2dt  \\
\ns & \ds + 2i\sum_{j=1}^n (\ell_{jt} +
\ell_{tj})(\ov_j v - v_j\ov)dt +
i(\Psi + \D \ell)(\ov dv - v d\ov)  \\
\ns & \ds  + \ell_t dv d\ov +
i\sum_{j=1}^n \ell_j (d\ov_j dv - dv_j
d\ov).
\end{array}
\end{equation}
Here
\begin{equation}\label{Id2eq1.1}
\begin{array}{ll}
M \3n& \ds = \b^2\ell_t |v|^2 +
i\b\sum_{j,k=1}^nb^{jk}\ell_j(\ov_kv-v_k\ov)\\
\ns & \ds = \ell_t |v|^2 +
i\sum_{j=1}^n \ell_j (\ov_j v - v_j
\ov);
\end{array}
\end{equation}
\begin{equation}\label{Id2eq1.2}
\begin{array}{ll}
A \3n& \ds
=\sum_{j,k=1}^nb^{jk}\ell_j\ell_k -
\sum_{j,k=1}^n(b^{jk}\ell_j)_k -\Psi \\
\ns & \ds  = \sum_{j=1}^n (\ell_j^2 -
\ell_{jj}) -\Psi;
\end{array}
\end{equation}
\begin{equation}\label{Id2eq1.3}
\begin{array}{ll}
D \3n & \ds =(\b^2\ell_t)_t
+\sum_{j,k=1}^n(b^{jk}\Psi_k)_j
+ 2\[\sum_{j,k=1}^n(b^{jk}\ell_j A)_k + A\Psi\]\\
\ns & \ds = \ell_{tt} + \sum_{j=1}^n
\Psi_{jj} + 2\sum_{j=1}^n (\ell_j A)_j
+ 2 A\Psi;
\end{array}
\end{equation}
\begin{equation}\label{Id2eq1.4}
\begin{array}{ll}
c^{jk} \3n &\ds  =
\sum_{j',k'=1}^n\[2(b^{j'k}\ell_{j'})_{k'}b^{jk'}
- (b^{jk}b^{j'k'}\ell_{j'})_{k'}\Psi\] - b^{jk}\\
\ns & \ds  =
\[2(b^{kk}\ell_{k})_{j}b^{jj} -
\sum_{j'=1}^n
(b^{jk}b^{j'j'}\ell_{j'})_{j'} -
b^{jk}\Psi\]\\
\ns & \ds = 2\ell_{jk} - \d^{jk}\D \ell
- \d^{jk}\Psi;
\end{array}
\end{equation}
and
\begin{equation}\label{Id2eq1.5}
\begin{array}{ll}
V_k \3n &\ds = -i
\b\sum_{j=1}^n\[b^{jk}\ell_j(vd\ov -\ov
 dv ) + b^{jk}\ell_t(v_j\ov-\ov_jv) dt\]\\
\ns & \ds \q - \Psi\sum_{j=1}^n b^{jk}(v_j\ov+\ov_jv) dt + \sum_{j=1}^n b^{jk}(2A\ell_j+\Psi_j)|v|^2 dt \\
\ns & \ds \q
+\sum_{j,j',k'=1}^n\(2b^{jk'}b^{j'k}-b^{jk}b^{j'k'}\)\ell_j(v_{j'}\ov_{k'}+\ov_{j'}v_{k'})
dt\\
\ns & \ds = -i\big[ \ell_k(vd\ov - \ov
dv) + \ell_t(v_j\ov -\ov_j v)dt \big] - \Psi(v_k\ov + \ov_k v)dt + (2A\ell_k + \Psi_k)|v|^2dt\\
\ns &  \ds \q + 2\sum_{j=1}^n \ell_j
(\ov_j v_k + v_j \ov_k)dt -
2\sum_{j'=1}^n \ell_k(v_j\ov_j)dt.
\end{array}
\end{equation}

\textbf{Step 2.} In this step, we
estimate the terms in the right-hand
side of the equality
\eqref{identity2.1} one by one.

First, from the definition of $\ell$,
$\f$(see \eqref{lvarphi}) and the
choice of $\psi$(see \eqref{psi}), we
have
\begin{equation}\label{lt1}
\begin{array}{ll}\ds
|\ell_t| & \ds = \Big| s\frac{2(2t-T)}{t^3(T-t)^3}\big( e^{4\l\psi} - e^{5\l |\psi|_{L^\infty(G)}} \big)  \Big| \\
\ns& \ds \leq \Big| s\frac{2(2t-T)}{t^3(T-t)^3} e^{5\l |\psi|_{L^\infty(G)}} \Big|  \\
\ns &\ds \leq  \Big| s\frac{C}{t^3(T-t)^3} e^{5\l \psi} \Big|\\
\ns & \ds  \leq
Cs\varphi^{1+\frac{1}{2}}
\end{array}
\end{equation}
and
\begin{equation}\label{ltt1}
\begin{array}{ll}
\ds |\ell_{tt}| & \ds  = \Big|  s\frac{20t^2 - 20tT + 6T^2}{t^4(T-t)^4} \big( e^{4\l\psi} - e^{5\l |\psi|_{L^\infty(G)}} \big) \Big| \\
\ns& \ds \leq \Big|  s\frac{C}{t^4(T-t)^4}  e^{5\l |\psi|_{L^\infty(G)}}  \Big| \\
\ns& \ds \leq \Big|  s\frac{C}{t^4(T-t)^4}  e^{8\l  \psi }  \Big|\\
\ns  &\ds  \leq Cs\f^2\leq Cs\f^3.
\end{array}
\end{equation}

We choose  below $\Psi = -\D \ell$,
then we have that
\begin{equation}\label{Id2eq2}
 A = \sum_{j=1}^n \ell_j^2 =  \sum_{j=1}^n \big(4s\l\f \psi_j )^2 =16s^2\l^2\varphi^2 |\nabla\psi|^2.
\end{equation}
Hence, we find
\begin{equation}\label{B}
\begin{array}{ll}\ds
D \3n & \ds = \ell_{tt} + \sum_{j=1}^n
\Psi_{jj} + 2\sum_{j=1}^n (\ell_j
A)_j + 2 A\Psi \\
\ns & \ds  = \ell_{tt} + \D(\D\ell) + 2\sum_{j=1}^n\big(4s\l\f\psi_j 16s^2\l^2\f^2|\nabla\psi|^2\big)_j - 32s^2\l^2\f^2|\nabla\psi|^2\D \ell  \\
\ns & \ds  =
384s^3\l^4\varphi^3|\nabla\psi|^4 -
\l^4\varphi O(s) - s^3\varphi^3 O(\l^3)
+ \ell_{tt}.
\end{array}
\end{equation}
Recalling that $x_0\in
(\mathbb{R}^n\setminus \overline G)$,
we know that
$$|\nabla\psi|>0\;\;\mbox{ in }\overline G.$$ From   \eqref{B} and  \eqref{ltt1}, we
know that there exists a $\l_0>0$ such
that for all $\l>\l_0$, one can find a
constant $s_0 = s_0(\l_0)$ so that for
any $s>s_0$, it holds that
\begin{equation}\label{B1}
D|v|^2 \geq
s^3\l^4\varphi^3|\nabla\psi|^4|v|^2.
\end{equation}
Since
$$
\begin{array}{ll}\ds
 c^{jk} = 2\ell_{jk} - \d^{jk}\D \ell - \d^{jk}\Psi  \\
\ns\ds\q\,\,\, = 32s\l^2\varphi\psi_j
\psi_k + 16s\l\varphi\psi_{jk},
\end{array}
$$
we see that
\begin{equation}\label{cjk}
\begin{array}{ll}\ds
\q \ds\sum_{j,k=1}^n c^{jk}(v_j\ov_k + v_k\ov_j)\\
\ns \ds  =
32s\l^2\varphi\sum_{j,k=1}^n\psi_j
\psi_k(v_j\ov_k + v_k\ov_j) +
16s\l\varphi
\sum_{j,k=1}^n \psi_{jk}(v_j\ov_k + v_k\ov_j)\\
\ns\ds =
32s\l^2\varphi\[\sum_{j=1}^n(\psi_jv_j)\sum_{k=1}^n
(\psi_k \ov_k) +
\sum_{k=1}^n(\psi_kv_k)\sum_{j=1}^n
(\psi_j \ov_j)  \] + 32s\l\varphi
\sum_{j=1}^n(v_j\ov_j
+ \ov_j v_j)\\
\ns   \ds  = 64s\l^2\varphi |\nabla\psi\cd\nabla v|^2 + 64 s\l\f |\nabla v|^2\\
\ns   \ds \geq  64 s\l\f |\nabla v|^2.
\end{array}
\end{equation}

Now we estimate the other terms in the
right-hand side of the equality
\eqref{identity2.1}. The first one
reads
\begin{equation}\label{ltj}
\begin{array}
{ll} \ds  2i\sum_{j=1}^n (\ell_{jt} +
\ell_{tj})(\ov_j v - v_j\ov) \3n & \ds
=
4i\sum_{j=1}^n s\l\psi_j \ell_t(\ov_j v - \ov v_j)\\
\ns & \ds \leq 2 s\varphi |\nabla v|^2
+ 2 s\l^2\varphi^3 |\nabla\psi|^2|v^2|.
\end{array}
\end{equation}
The second one satisfies
\begin{eqnarray}\label{liiPsi}
i(\Psi + \D \ell)(\ov dv - v d\ov) = 0.
\end{eqnarray}

For  estimating  the third and the
fourth one, we need to take mean value
and get that
\begin{equation}\label{dvdov}
\begin{array}{ll}\ds
\mathbb{E}\big(\ell_t dv d\ov\big) \3n&
\ds= \mathbb{E}\big[\ell_t(\theta
\ell_t zdt + \theta
dz)(\overline{\theta \ell_t zdt +
\theta dz)}\big] = \mathbb{E}(\ell_t
\theta^2 dz d\bar{z})
\\
\ns & \ds \leq 2s\theta^2
\varphi^{\frac{3}{2}}\mathbb{E}|Z|^2dt.
\end{array}
\end{equation}
Here we utilize the inequality \eqref{lt1}.

\vspace{0.1cm}

Further,
$$
\begin{array}{ll}\ds
\mathbb{E}(d\ov_j dv) \3n &
= \mathbb{E}\big[\overline{\big( \theta \ell_t z dt + \theta dz  \big)}_j \big( \theta \ell_t z dt + \theta dz  \big)\big] \\
\ns& \ds  = \mathbb{E} \big[\, \overline{(\theta dz)}_j (\theta dz) \big]\\
\ns& \ds =  \mathbb{E} \big[\, \overline{\big( s\l\f\psi_j\theta dz + \theta dz_j   \big)}\theta dz  \big]\\
\ns & \ds  = s\l\f\psi_j\theta^2 \mathbb{E}d\bar z dz + \theta^2 \mathbb{E}d\bar z_j dz \\
\ns & \ds  =  s\l\f\psi_j\theta^2
\mathbb{E}|Z|^2dt + \theta^2
\mathbb{E}\big(\overline{Z}_j Z\big)dt.
\end{array}
$$
Similarly, we can get that
$$
\begin{array}{ll}\ds
\mathbb{E}(dv_j
d\ov)=s\l\f\psi_j\theta^2
\mathbb{E}|Z|^2dt + \theta^2
\mathbb{E}\big(\overline{Z} Z_j\big)dt.
\end{array}
$$
Therefore, the fourth one satisfies
that
\begin{equation}\label{dvjdv}
\begin{array}{ll}
\ds \q i\mathbb{E}\sum_{j=1}^n \ell_j
(d\ov_j dv - dv_j d\ov)=
is\l\f\psi_j\big[\theta^2
\mathbb{E}\big(\overline{Z}_j Z\big)dt
-
 \theta^2
\mathbb{E}\big(\overline{Z}_j
Z\big)dt\big].
\end{array}
\end{equation}

\textbf{Step 3.} Integrating the
equality \eqref{identity2.1} in $Q$,
taking mean value in both sides, and
noting \eqref{Id2eq2}--\eqref{dvjdv},
we obtain that
\begin{equation}\label{inep1}
\begin{array}{ll}
\ds \q\mathbb{E}\int_Q
\Big(s^3\l^4\varphi^3 |v|^2 +
s\l^2\varphi |\nabla v|^2\Big) dxdt +
2\mathbb{E}\int_Q \Big|\!- i
\ell_t v - 2\sum_{j=1}^n  \ell_j v_j  + \Psi v\Big|^2dxdt\\
\ns \ds  \leq  \mathbb{E}\int_Q \[
\theta\cP y {\Big( i \ell_t \bar{v} -
2\sum_{j =1}^n  \ell_j \bar{v}_j + \Psi
\bar{v}\Big)} + \theta\overline{\cP y}
{\Big(- i \ell_t v - 2\sum_{j =1}^n
\ell_j v_j + \Psi v\Big)} \]dx\\
\ns \ds \q  +\; C\mathbb{E}\int_Q
\theta^2\big(s^2\l^2 \varphi^2|Z| +
 |\nabla Z|^2\big)dxdt + \mathbb{E}\int_Q dM dx +
\mathbb{E}\int_Q \div V dx.
\end{array}
\end{equation}

Now we analyze the  terms in the
right-hand side of the inequality
\eqref{inep1} one by one.

The first one reads
\begin{equation}\label{intprin}
\begin{array}{ll} \ds
\q\mathbb{E}\int_Q \[ \theta\cP z
{\Big( i \ell_t \bar{v} - 2\sum_{j
=1}^n \ell_j \bar{v}_j + \Psi
\bar{v}\Big)} + \theta\overline{\cP z}
{\Big(- i \ell_t v -
2\sum_{j =1}^n  \ell_j v_j + \Psi v\Big)} \]dx \\
\ns  \ds  =  \ds \mathbb{E}\int_Q \[
\theta \big(b_1 \cdot \nabla z + b_2 z
+ b_3Z \big) {\( i \ell_t \bar{v} -
2\sum_{j =1}^n \ell_j
\bar{v}_j + \Psi \bar{v}\)}\\
\ns\ds \qq\qq + \theta {\big(-b_1 \cdot
\nabla \bar z + b_2 \bar z +
b_3\overline Z \big)} {\Big(- i \ell_t
v - 2\sum_{j =1}^n  \ell_j v_j  + \Psi
v\Big)} \] dxdt \\
\ns \ds \leq  2\mathbb{E}\int_Q
\Big\{\theta^2\big|b_1 \cdot \nabla z +
b_2 z + b_3Z \big|^2 + \Big|- i\b
\ell_t v - 2\sum_{j =1}^n \ell_j v_j +
\Psi v\Big|^2 \Big\}dxdt.
\end{array}
\end{equation}

From the choice of $\theta$, we know
that $v(0)=v(T)=0$. Hence, we have
\begin{equation}\label{idm}
\int_Q dM dx = 0.
\end{equation}
Further, by Stokes' Theorem  and noting
that $v=z=0$ on $(0,T)\times\G$, we
find
\begin{equation}\label{intV}
\begin{array}
{ll} \ds \mathbb{E}\int_Q \div V dx\3n
&\ds = \ds \mathbb{E}\int_{\Si}
2\sum_{k=1}^n\sum_{j=1}^n\big[
 \ell_j(\ov_j v_k + v_j \ov_k)\nu^k - \ell_k \nu_k v_j \ov_j
\big]d\Si\\
\ns  &\ds =  \ds \mathbb{E}\int_{\Si}
\Big(4\sum_{j=1}^n \ell_j \nu_j \Big|
\frac{\pa v}{\pa \nu} \Big|^2 -
2\sum_{k=1}^n \ell_k \nu_k \Big|
\frac{\pa v}{\pa \nu} \Big|^2\Big) d\Si\\
\ns  &= \ds  \mathbb{E}\int_{\Si}
2\sum_{k=1}^n \ell_k \nu_k \Big|
\frac{\pa v}{\pa \nu} \Big|^2 d\Si \\
\ns &\ds \leq  C\mathbb{E}\int_0^T
\int_{\G_0} \theta^2 s\l\varphi \Big|
\frac{\pa z}{\pa \nu} \Big|^2 d\G dt.
\end{array}
\end{equation}
From \eqref{inep1}--\eqref{intV}, we have
\begin{equation}\label{car1}
\begin{array}{ll}
\q \ds \mathbb{E}\int_Q
\Big(s^3\l^4\varphi^3 |v|^2 +
s\l\varphi
|\nabla v|^2\Big) dxdt \\
\ns   \ds \leq C\,\mathbb{E}\int_Q
\theta^2 |b_1 \cdot \nabla z + b_2 z +
b_3Z|^2 dxdt +
C\,\mathbb{E}\int_0^T\int_{\G_0}\theta^2
s\l\varphi\Big| \frac{\pa z}{\pa \nu}\Big|^2d\G dt\\
\ns \ds   \q  +\, C\mathbb{E}\int_Q
\theta^2\big(s^2\l^2 \varphi^2|Z|^2 +
|\nabla Z|^2\big) dxdt.
\end{array}
\end{equation}
Noting that $z_i = \theta^{-1}(v_i -
\ell_i v) = \theta^{-1}(v_i -
s\l\varphi\psi_i v)$, we get
\begin{equation}\label{vtoy}
\theta^2\big(|\nabla z|^2 +
s^2\l^2\varphi^2 |z|^2\big)\leq
C\big(|\nabla v|^2 + s^2\l^2\varphi^2
|v|^2\big).
\end{equation}
Therefore, it follows from \eqref{car1}
that
\begin{equation}\label{car2}
\begin{array}{ll}
\ds \q\mathbb{E}\int_Q
\big(s^3\l^4\varphi^3 |z|^2 +
s\l\varphi
|\nabla z|^2\big) dxdt \\
\ns \ds   \leq  C\mathbb{E}\int_Q
\theta^2\big( |b_1|^2 | \nabla z|^2 +
b_2^2 |z|^2 + b_3^2|Z|^2\big) dxdt +
C\mathbb{E}\int_0^T\int_{\G_0}\theta^2
s\l\varphi\Big| \frac{\pa
z}{\pa \nu}\Big|^2d\G dt \\
\ns  \ds \q + C\mathbb{E}\int_Q
\theta^2\big(s^2\l^2 \varphi^2|Z|^2 +
|\nabla Z|^2\big)  dxdt.
\end{array}
\end{equation}

Taking   $\l_1 =\l_0$ and $s_1 =
\max(s_0, Cr_1)$, and utilizing the
inequality \eqref{car2}, we conclude
the desired inequality \eqref{carleman
est}.
\endpf

\section{Proof of Theorem \ref{observability}}

In this section, we prove Theorem
\ref{observability}  by virtue  of
Theorem \ref{thcarleman est}.

\vspace{0.1cm}

 {\it Proof of
Theorem \ref{observability}}: By means
of the definition of $\ell$ and
$\theta$(see \eqref{lvarphi}), it holds
that
\begin{equation}\label{final1}
\begin{array}
{ll} &\ds \mathbb{E}\int_Q
\theta^2\Big(\varphi^3 |z|^2 + \varphi
|\nabla z|^2\Big) dxdt\\
\ns  \geq & \ds
\min_{x\in\overline{G}}\Big(\varphi\Big(\frac{T}{2},x\Big)
\theta^2\Big(\frac{T}{4},x\Big)\Big)\mathbb{E}\int_{\frac{T}{4}}^{\frac{3T}{4}}\int_G\big(|z|^2+|\nabla
z|^2\big)dxdt,
\end{array}
\end{equation}

\begin{equation}\label{final2}
\begin{array}
{ll} &\ds \mathbb{E}\int_Q \theta^2(
\varphi^2|Z|^2 + |\nabla
Z|^2)dxdt \\
\ns \leq & \ds \max_{(x,t)\in
\overline{Q}}\big(\varphi^2(t,x)\theta^2(t,x)\big)
\mathbb{E}\int_Q\big( |Z|^2 + |\nabla
Z|^2\big)dxdt
\end{array}
\end{equation}
and that
\begin{equation}\label{final3}
\mathbb{E}\int_0^T\int_{\G_0}\theta^2
\varphi\Big| \frac{\pa z}{\pa
\nu}\Big|^2d\G dt \leq \max_{(x,t)\in
\overline{Q}}\big(\varphi(t,x)\theta^2(t,x)\big)\mathbb{E}\int_0^T\int_{\G_0}
\Big| \frac{\pa z}{\pa \nu}\Big|^2d\G
dt.
\end{equation}

From  \eqref{carleman est} and
\eqref{final1}--(\ref{final3}), we
deduce that
\begin{eqnarray}\label{final4}
\begin{array}
{ll} &\ds
\mathbb{E}\int_{\frac{T}{4}}^{\frac{3T}{4}}\int_G(|z|^2+|\nabla
z|^2)dxdt\\
\ns \leq & \ds  C r_1
\frac{\max_{(x,t)\in
\overline{Q}}\Big(\varphi^2(t,x)\theta^2(t,x)\Big)}{\min_{x\in\overline{G}}\Big(\varphi(\frac{T}{2},x)\theta^2(\frac{T}{4},x)\Big)}\\
\ns & \ds \q\times\left\{
\mathbb{E}\int_Q(|Z|^2 + |\nabla
Z|^2)dxdt +
\mathbb{E}\int_0^T\int_{\G_0} \Big|
\frac{\pa z}{\pa
\nu}\Big|^2d\G dt\right\}\\
\ns \leq & \ds   e^{ Cr_1 }\left\{
\mathbb{E}\int_Q(|Z|^2 + |\nabla
Z|^2)dxdt +
\mathbb{E}\int_0^T\int_{\G_0} \Big|
\frac{\pa z}{\pa \nu}\Big|^2d\G
dt\right\}.
\end{array}
\end{eqnarray}

Utilizing (\ref{final4}) and
(\ref{energyesti1}), we obtain that
\begin{eqnarray}\label{final5}
\begin{array}
{ll} &\ds \mathbb{E}\int_G(|z_T|^2 + |\nabla z_T|^2)dx \\
\ns \leq & \ds   e^{C r_1 }\left\{
\mathbb{E}\int_Q( |Z|^2 + |\nabla
Z|^2)dxdt +
\mathbb{E}\int_0^T\int_{\G_0} \Big|
\frac{\pa z}{\pa \nu}\Big|^2d\G
dt\right\},
\end{array}
\end{eqnarray}
which concludes Theorem
\ref{observability} immediately.\endpf

\section{Proof of Theorem \ref{exact th}--\ref{non control
th}}

This section is addressed to  proofs of
Theorem \ref{exact th}-\ref{non control
th}.

\vspace{0.2cm}

{\it Proof of Theorem \ref{exact
th}}\,: Since the system
\eqref{system1} is linear, we only need
to show that the attainable set $A_T$
at time $T$ with initial datum $y(0)=0$
is $L^2_{\cF_T}(\O;H^{-1}(G))$, that
is, for any $y_1\in
L^2_{\cF_T}(\O;H^{-1}(G))$, we can find
a pair of control
$$(u,g)\in
L_{\cF}^2(0,T; L^2(\G_{0}))\times
L_{\cF}^2(0,T;H^{-1}(G))$$ such that
the solution to the system
\eqref{system1} with $y(0)=0$ satisfies
that $y(T)=y_1$ in
$L^2_{\cF_T}(\O;H^{-1}(G))$. We achieve
this goal by duality argument.

Let
$$
b_1 = -a_1,\q b_2 = -\div(a_1)+a_2,
b_3=-a_3
$$ in the equation
\eqref{system2}. We set
$$
\begin{array}{ll}\ds
\cX\=\Big\{\(\frac{\pa
z}{\pa\nu}\Big|_{\G_{0}},
Z\)\;\Big|\;(z,Z)\hb{ solves the
equation
 }\eqref{system2}\mbox{ with some } \\
\ns\ds \hspace{3.9cm} z_T\in
L^2(\O,\cF_T,P;H_0^1(G))\Big\}.
\end{array}
$$
Clearly, $\cX$ is a linear subspace of
$L_{\cF}^2(0,T;L^2(\G_{0}))\times
L_{\cF}^2(0,T;H^{1}_0(G))$. Let us
define a linear functional $\cL$ on
$\cX$ as follows:
$$
\cL\(\frac{\pa
z}{\pa\nu}\Big|_{\G_{0}},
Z\)=\mathbb{E}\langle
y_1,z_T\rangle_{H^{-1}(G),H_0^1(G)}  -
\mE\int_0^T\int_G z f dxdt.
$$
From Theorem \ref{observability}, we
see that $\cL$ is a bounded linear
functional on $\cX$. By means of the
Hahn-Banach theorem, $\cL$ can be
extended to be a bounded linear
functional on the space $L_{\cF}^2(0,T;
L^2(\G_{0}))\times
L_{\cF}^2(0,T;H^{-1}(G))$. For
simplicity, we still use $\cL$ to
denote this extension. Now, by  Riesz
Representation theorem, we know that
there is a pair of random fields
$$(u,
g)\in L_{\cF}^2(0,T; L^2(\G_{0}))\times
L_{\cF}^2(0,T;H^{-1}(G))
$$
such that
\begin{equation}\label{con eq0}
\begin{array}{ll}\ds
\q\mE\langle
y_1,z_T\rangle_{H^{-1}(G),H_0^1(G)}  -
\mE\int_0^T\int_G  zf dxdt \\
\ns\ds = \mE\int_0^T\int_{\G_{0}}
\frac{\pa z}{\pa\nu}u d\G dt +
\mE\int_0^T \langle g, Z
\rangle_{H^{-1}(G),H_0^1(G)} dt.
\end{array}
\end{equation}

We claim that this pair of random
fields $(u,g)$ is  the control we need.
In fact, from the definition of the
solution to \eqref{system1}, we have
\begin{equation}\label{con eq1}
\begin{array}{ll}\ds
\q\mE\langle
y(T),z_T\rangle_{H^{-1}(G),H_0^1(G)}\\
\ns\ds
 = \mE\int_0^T\int_G  zf dxdt
 + \mE\int_0^T\int_{\G_{0}} \frac{\pa
z}{\pa\nu}u d\G dt + \mE\int_0^T
\langle g,
Z\rangle_{H^{-1}(G),H_0^1(G)} dt.
\end{array}
\end{equation}
From   \eqref{con eq0} and  \eqref{con
eq1}, we see
\begin{equation}\label{con eq3}
\mE\langle y_1,
z_T\rangle_{H^{-1}(G),H_0^1(G)} =
\mE\langle y(T,\cd),z_T
\rangle_{H^{-1}(G),H_0^1(G)}.
\end{equation}
Since $z_T$ can be arbitrary element in
$L^2_{\cF_T}(\O;H_0^1(G))$, from the
equality \eqref{con eq3}, we get
$y(T)=y_1$ in $H^{-1}(G)$, $P$-a.s.
\endpf

At last, we prove Theorem \ref{non
control th}. In order to present the
key idea in the simplest way, we only
consider a very special case of the
system \eqref{system1}, that is,
$a_1=0$, $a_2=0$ and $a_3 = 1$. The
argument for the general case is very
similar.

\vspace{0.2cm}

{\it Proof of Theorem \ref{non control
th}}\,: Let us assume that $u\equiv 0$.
In this case, the system
\eqref{system1} is
\begin{equation}\label{system1.1}
\left\{
\begin{array}{ll}\ds
idy + \D ydtdt =   (y+g)dB(t)   &\mbox{ in } (0,T)\t G,\\
\ns\ds y(t,0) = 0 & \mbox{ on } (0,T)\t\G_0,\\
\ns\ds y(0) = y_0 &\mbox{ in } G.
\end{array}
\right.
\end{equation}
Since the system \eqref{system1.1} is
linear, we only need to show that the
attainable set $A_T$ at time $T$ for
the initial datum $y_0 =0$ is not
$L^2_{\cF_T}(\O;H^{-1}(G))$. The
solution of the system
\eqref{system1.1} is
\begin{equation}\label{sol2}
y(T) = S(T)y_0 - i\int_0^T
S(T-s)\big[y(s)+g(s)\big]dB(s).
\end{equation}
Here $\{S(t)\}_{t\geq 0}$ is the
semigroup generated by the following
operator
$$
\left\{
\begin{array}{ll}\ds
D(A)=H^2(G)\cap H_0^1(G),\\
\ns\ds A\f = i\D\f,\q\forall\f\in D(A).
\end{array}
\right.
$$
From  \eqref{sol2}, we find $\mE (y(T))
= \mE(S(T)y_0)$. Thus, if we choose a
$y_1 \in L^2_{\cF_T}(\O;H^{-1}(G))$
such that $\mE (y_1) \neq 0$, then
$y_1$ is not in $A_T$, which completes
the proof.
\endpf

\section{Further comments and open
problems}\label{comment}

There are plenty of open problems in
the topic of this paper. Some of them
are particularly relevant and could
need important new ideas and further
developments:

\begin{itemize}

\item {\bf Null and approximate controllability for
stochastic Schr\"{o}ding equations}

In this paper, we study the exact
controllability for stochastic
Schr\"{o}dinger equations. As immediate
consequences, we can obtain the null
and approximate controllability for the
same system. However, in order to get
these two kinds of controllability, we
have no reasons to use two controls. By
the proof of Theorem \ref{non control
th}, we know that it is not enough to
put one control in the diffusion term
to get the null or approximate
controllability. On the other hand,
suggested by the result in \cite{Lu0},
we believe one boundary control in the
drift term can guarantee the null and
approximate controllability of
\eqref{system1}. If we want to prove
this by following the method in this
paper, we will meet some essential
difficulty. For example, to get the
null controllability, we should prove
the following inequality for the
solution to \eqref{system2}
\begin{equation}\label{4.16-eq4}
\mE|z(0)|^2_{H_0^1(G)} \leq
C\int_0^T\int_{\G_0} \Big| \frac{\pa
z}{\pa\nu} \Big|^2 d\G dt.
\end{equation}
However, if we utilize the method in
this paper, we only get
$$
\mE|z(0)|^2_{H_0^1(G)} \leq
C\int_0^T\int_{\G_0} \Big| \frac{\pa
z}{\pa\nu} \Big|^2 d\G dt + \int_0^T
|Z|_{H_0^1(G)}^2 dt.
$$
There is an additional term containing
$Z$ in the right hand side. This terms
comes from the fact that, in the
Carleman estimate, we regard $Z$ as an
nonhomogeneous term rather than part of
the solution. Hence, it has to appear
in the right hand side of the
inequality. Therefore, we believe that
one should introduce some new
technique, for example, a Carleman
estimate in which the fact $Z$ is part
of the solution is involved, to get rid
of the additional term containing $Z$.
However, we do not know how to achieve
this goal now.

\item{\bf Exact controllability for stochastic Schr\"{o}dinger equations
with less restrictive condition}

In this paper, we get the exact
controllability for \eqref{system1} for
$\G_0$ given by \eqref{G0}. It is well
known that a sharp sufficient condition
for exact controllability for
deterministic Schr\"{o}dinger equations
in analytic domain with time invariant
lower order terms is that the triple
$(G, \G_0, T)$ satisfies the Geometric
Optic Condition introduced in
\cite{Bardos-Lebeau-Rauch1}(see
\cite{Lebeau} for example). It would be
quite interesting and challenging to
extend this result to the stochastic
setting, but it seems that there are
lots of things should be done before
solving this problem. For instance, the
main idea in \cite{Lebeau} is as
follows:

Based on a diadic decomposition of the
Fourier representation of solutions of
the Schr\"{o}dinger equation, the
author show that, to some extent, one
can regards them as superposition of an
infinite sequence of solutions of wave
equations with velocity of propagation
tending to infinity. Then, the fact
that the Geometric Optic Condition is
satisfied for some finite time $T$
suffices for the exact controllability
of the Schr\"{o}dinger equation  to
hold for all $T>0$.

If one follow the above idea, then
propagation of singularities for
stochastic partial differential
equations, at least, for stochastic
hyperbolic equations, should be
established. However, as far as we
know, this topic is completely open.

Further, there are some results showing
that, in some situations in which the
Geometric Optic Condition is not
fulfilled in any time $T$, one can
still obtain exact controllability for
the Schr\"{o}dinger equation. For
instance, in \cite{Jaffard}, it is
showed that, when the domain $G$ is a
square, for any open non-empty subset
$G_0$ of $G$, the exact controllability
of the Schr\"{o}dinger equation  holds
in any time $T$, in the space $L^2(G)$
and with internal controls in $L^2
((0,T)\times G_0)$. How to prove such
kind of result is very interesting.

\end{itemize}



\end{document}